\documentclass[11pt]{amsart}
\usepackage{amsopn}
\usepackage{amssymb, amscd}

\newcommand{\nc}{\newcommand}

\nc{\fg}{\mathfrak{f} } \nc{\vg}{\mathfrak{v} } \nc{\wg}{\mathfrak{w} }
\nc{\zg}{\mathfrak{z} } \nc{\ngo}{\mathfrak{n} } \nc{\kg}{\mathfrak{k} }
\nc{\mg}{\mathfrak{m} } \nc{\bg}{\mathfrak{b} } \nc{\ggo}{\mathfrak{g} }
\nc{\ggob}{\overline{\mathfrak{g}} } \nc{\sog}{\mathfrak{so} }
\nc{\sug}{\mathfrak{su} } \nc{\spg}{\mathfrak{sp} } \nc{\slg}{\mathfrak{sl} }
\nc{\glg}{\mathfrak{gl} } \nc{\cg}{\mathfrak{c} } \nc{\rg}{\mathfrak{r} }
\nc{\hg}{\mathfrak{h} } \nc{\tg}{\mathfrak{t} } \nc{\ug}{\mathfrak{u} }
\nc{\dg}{\mathfrak{d} } \nc{\ag}{\mathfrak{a} } \nc{\pg}{\mathfrak{p} }
\nc{\sg}{\mathfrak{s} } \nc{\affg}{\mathfrak{aff} }

\nc{\pca}{\mathcal{P}} \nc{\nca}{\mathcal{N}} \nc{\lca}{\mathcal{L}}
\nc{\oca}{\mathcal{O}} \nc{\mca}{\mathcal{M}} \nc{\tca}{\mathcal{T}}
\nc{\aca}{\mathcal{A}} \nc{\cca}{\mathcal{C}} \nc{\gca}{\mathcal{G}}
\nc{\sca}{\mathcal{S}} \nc{\hca}{\mathcal{H}} \nc{\bca}{\mathcal{B}}
\nc{\dca}{\mathcal{D}} \nc{\val}{\operatorname{val}}

\nc{\vp}{\varphi} \nc{\ddt}{\tfrac{{\rm d}}{{\rm d}t}}
\nc{\dpar}{\tfrac{\partial}{\partial t}} \nc{\im}{\mathtt{i}}

\nc{\SO}{\mathrm{SO}} \nc{\Spe}{\mathrm{Sp}} \nc{\Sl}{\mathrm{SL}}
\nc{\SU}{\mathrm{SU}} \nc{\Or}{\mathrm{O}} \nc{\U}{\mathrm{U}} \nc{\Gl}{\mathrm{GL}}
\nc{\Se}{\mathrm{S}} \nc{\Cl}{\mathrm{Cl}} \nc{\Spein}{\mathrm{Spin}}
\nc{\Pin}{\mathrm{Pin}} \nc{\G}{\mathrm{GL}_n(\RR)} \nc{\g}{\mathfrak{gl}_n(\RR)}

\nc{\RR}{{\Bbb R}} \nc{\HH}{{\Bbb H}} \nc{\CC}{{\Bbb C}} \nc{\ZZ}{{\Bbb Z}}
\nc{\FF}{{\Bbb F}} \nc{\NN}{{\Bbb N}} \nc{\QQ}{{\Bbb Q}} \nc{\PP}{{\Bbb P}}

\nc{\vs}{\vspace{.2cm}} \nc{\vsp}{\vspace{1cm}} \nc{\ip}{\langle\cdot,\cdot\rangle}
\nc{\ipp}{(\cdot,\cdot)} \nc{\la}{\langle} \nc{\ra}{\rangle} \nc{\unm}{\tfrac{1}{2}}
\nc{\unc}{\tfrac{1}{4}} \nc{\und}{\tfrac{1}{16}} \nc{\no}{\vs\noindent}
\nc{\lamn}{\Lambda^2(\RR^n)^*\otimes\RR^n} \nc{\lamp}{\Lambda^2\pg^*\otimes\pg}
\nc{\lamg}{\Lambda^2\ggo^*\otimes\ggo} \nc{\lamngo}{\Lambda^2\ngo^*\otimes\ngo}
\nc{\tangz}{{\rm T}^{\rm Zar}} \nc{\mum}{/\!\!/} \nc{\kir}{/\!\!/\!\!/}
\nc{\Ri}{\tfrac{4\Ric_{\mu}}{||\mu||^2}} \nc{\ds}{\displaystyle}
\nc{\ben}{\begin{enumerate}} \nc{\een}{\end{enumerate}} \nc{\f}{\frac}
\nc{\lb}{[\cdot,\cdot]} \nc{\isn}{\tfrac{1}{||v||^2}}
\nc{\gkp}{(\ggo=\kg\oplus\pg,\ip)} \nc{\ukh}{(\ug=\kg\oplus\hg,\ip)}

\nc{\Hess}{\operatorname{Hess}} \nc{\ad}{\operatorname{ad}}
\nc{\Ad}{\operatorname{Ad}} \nc{\rank}{\operatorname{rank}}
\nc{\Irr}{\operatorname{Irr}} \nc{\End}{\operatorname{End}}
\nc{\Aut}{\operatorname{Aut}} \nc{\Inn}{\operatorname{Inn}}
\nc{\Der}{\operatorname{Der}} \nc{\Ker}{\operatorname{Ker}}
\nc{\Iso}{\operatorname{Iso}} \nc{\Diff}{\operatorname{Diff}}
\nc{\Lie}{\operatorname{Lie}} \nc{\tr}{\operatorname{tr}} \nc{\dif}{\operatorname{d}}
\nc{\sen}{\operatorname{sen}} \nc{\modu}{\operatorname{mod}}
\nc{\Riem}{\operatorname{Rm}} \nc{\Ricci}{\operatorname{Ric}}
\nc{\sym}{\operatorname{sym}} \nc{\symac}{\operatorname{sym^{ac}}}
\nc{\symc}{\operatorname{sym^{c}}} \nc{\scalar}{\operatorname{R}}
\nc{\grad}{\operatorname{grad}} \nc{\ricci}{\operatorname{Rc}}
\nc{\nr}{\operatorname{nr}} \nc{\riccic}{\operatorname{ric^{c}}}
\nc{\riccig}{\operatorname{ric^{\gamma}}} \nc{\Rin}{\operatorname{M}}
\nc{\Le}{\operatorname{L}} \nc{\tang}{\operatorname{T}}
\nc{\level}{\operatorname{level}} \nc{\rad}{\operatorname{r}}
\nc{\abel}{\operatorname{ab}} \nc{\CH}{\operatorname{CH}}
\nc{\mcc}{\operatorname{mcc}} \nc{\Adj}{\operatorname{Adj}}
\nc{\Order}{\operatorname{O}}

\theoremstyle{plain}
\newtheorem{theorem}{Theorem}[section]
\newtheorem{proposition}[theorem]{Proposition}
\newtheorem{corollary}[theorem]{Corollary}
\newtheorem{lemma}[theorem]{Lemma}

\theoremstyle{definition}
\newtheorem{definition}[theorem]{Definition}

\theoremstyle{remark}
\newtheorem{remark}[theorem]{Remark}

\newtheorem{example}[theorem]{Example}

\title{The Ricci flow for simply connected nilmanifolds}

\author{Jorge Lauret}

\thanks{This research was partially supported by grants from CONICET, ANPCyT (Argentina), SeCyT (UN C\'ordoba) and UC Berkeley.}

\address{FaMAF and CIEM, Universidad Nacional de C\'ordoba, C\'ordoba, Argentina}
\email{lauret@famaf.unc.edu.ar}

\begin{document}

\maketitle

\begin{abstract}
We prove that the Ricci flow $g(t)$ starting at any metric on $\RR^n$ that is invariant by a transitive nilpotent Lie group $N$ can be obtained by solving an ODE for a curve  of nilpotent Lie brackets on $\RR^n$.  By using that this ODE is the negative gradient flow of a homogeneous polynomial, we obtain that $g(t)$ is type-III, and, up to pull-back by time-dependent diffeomorphisms, that $g(t)$ converges to the flat metric, and the rescaling $|\!\scalar(g(t))|\,g(t)$ converges to a Ricci soliton in $C^\infty$, uniformly on compact sets in $\RR^n$.  The Ricci soliton limit is also invariant by some transitive nilpotent Lie group, though possibly nonisomorphic to $N$.
\end{abstract}


\section{Introduction}

It is expected that the Ricci flow behaves nicely in the presence of an additional structure on the starting metric.  This is well known in the case of K\"ahler manifolds, warped products, locally homogeneous $3$-manifolds and rotationally symmetric manifolds, among many others.  We study in this paper the case of nilmanifolds, a structure with strong algebraic data involved but yet geometrically very rich and exotic.

Let $g$ be a Riemannian metric on $\RR^n$ that is invariant under a transitive nilpotent Lie group $N$.  Let us consider the Ricci flow $g(t)$ starting at $g$, that is,
$$
\dpar g(t)=-2\ricci(g(t)),\qquad g(0)=g,
$$
and $g(t)$ is $N$-invariant for all $t$.  We prove the following results:

\begin{itemize}
\item $g(t)$ is defined for $t\in[0,\infty)$ and there exists a constant $C_n$ that only depends on $n$ such that
$\|\Riem(g(t))\|\leq \frac{C_n}{t}$ for all $t\in(0,\infty)$; in particular, $g(t)$ is a type-III solution.

\item The quantity $\frac{\|\ricci(g(t))\|}{|R(g(t))|}$ is strictly decreasing for all $t$ unless $g$ is a Ricci soliton.

\item If $\{ h(t)\}\subset\G\subset\Diff(\RR^n)$ is the solution to the ODE
$$
\ddt h(t)=-h(t)\Ricci_t,\qquad h(0)=I,
$$
where $\Ricci_t=\Ricci(g(t))(0)$ is the Ricci operator of $g(t)$ at the point $0\in\RR^n$, then the metrics $h(t)^*g(t)$ converge in $C^\infty$ to the flat metric uniformly on compact sets in $\RR^n$, as $t\to\infty$.

\item $h(t)^*g(t)$ is the negative gradient flow of the functional square norm of the Ricci tensor.

\item If we denote by $\scalar_t:=\scalar(g(t))$ the scalar curvature of $g(t)$ (recall it is constant on $\RR^n$), and $g$ is nonflat, then the metrics
    $$\tilde{g}(t):=|\!\scalar_t\!| h(t)^*g(t)$$
    converges in $C^\infty$ to a Ricci soliton metric $g_\infty$, uniformly on compact sets in $\RR^n$, as $t\to\infty$.  The metric $g_\infty$ on $\RR^n$ is also invariant under a transitive nilpotent Lie group, though possibly non-isomorphic to $N$.  We note that  $\scalar(\tilde{g}(t))\equiv -1$ and so $g_{\infty}$ is never flat as $R(g_\infty)=-1$.
\end{itemize}

These results show that most of the nice properties of the Ricci flow for nilmanifolds of dimension $3$ and $4$ discovered and proved in \cite{IsnJck,IsnJckLu,Glc,Ltt}, actually hold in the general case.

Our approach is based on an ODE for Lie brackets which is equivalent in a natural and specific sense to the Ricci flow $g(t)$ starting at any nilmanifold, and may be described as follows.

Let $\mca_n=(\mca_n,C^\infty)$ be the space of Riemannian metrics on $\RR^n$ endowed with the compact open $C^\infty$-topology (i.e. uniform $C^\infty$-convergence on compact sets ).  For each $\mu\in\nca_n\subset\lamn$, the algebraic subset of nilpotent Lie brackets on $\RR^n$, a metric $g_\mu\in\mca_n$ that is invariant by some transitive nilpotent Lie group can be defined in such a way that
$$
\nca_n\hookrightarrow\mca_n, \qquad\mu\mapsto g_\mu,
$$
turns out to be an embedding (i.e. $\mu_k\to\lambda$ (vector space topology) if and only if $g_{\mu_k}\to g_\lambda$ in $C^\infty$).  $\nca_n$ contains, up to isometry, all left-invariant metrics on simply connected nilpotent Lie groups of dimension $n$ (see Section \ref{nil}).

Let $g(t)$ be the Ricci flow with $g(0)=g_{\mu_0}$, $\mu_0\in\nca_n$.  We prove that $g(t)$ is given, up to pull-back by time-dependent diffeomorphisms, by $g_{\mu(t)}$ for the solution $\mu(t)\in\subset\nca_n$ to the ODE ({\it bracket flow})
$$
\ddt\mu(t)=\mu(t)(\Ricci_{\mu(t)}\cdot,\cdot) + \mu(t)(\cdot,\Ricci_{\mu(t)}\cdot) -\Ricci_{\mu(t)}\mu(t)(\cdot,\cdot), \qquad \mu(0)=\mu_0,
$$
where $\Ricci_\mu:\RR^n\longrightarrow\RR^n$ denotes the Ricci operator of $g_\mu$ at $0\in\RR^n$ (see Theorem \ref{eqfl}).  This ODE is precisely the negative gradient flow of the $4$-degree homogeneous polynomial $F(\mu)=\tr{\Ricci_\mu^2}$, whose only critical point is the global minima $0\in\nca_n$ (i.e. the flat metric $g_0$).  The fact that $g(t)$ is a type-III solution follows from an upper bound for the rate of decay of $\|\mu(t)\|\to 0$ (see Section \ref{t3}).

However, when restricted to any sphere, the critical points of $F$ are precisely the Ricci soliton metrics in $\nca_n$, and the convergence to a unique Ricci soliton of the corresponding normalized Ricci flow (the one for which $\scalar(g_\mu)=-\unc\|\mu\|^2\equiv -1$) follows from the known Thom conjecture (see Section \ref{scn}).

The bracket flow has already been used to study the Ricci flow for nilmanifolds in \cite{Gzh} and \cite{Pyn}, and for $3$-dimensional unimodular Lie groups in \cite{GlcPyn} (see Remark \ref{GP}). In \cite{homRF}, the equivalence (up to pull-back by time-dependent diffeomorphisms) between the Ricci and bracket flows has been proved for homogeneous Riemannian manifolds in general.

\vs \noindent {\it Acknowledgements.}  Part of this research was carried out while
visiting the University of California at Berkeley. I am very grateful to John Lott for the invitation and for fruitful discussions on the
topic of this paper.

\section{Nilmanifolds as metrics on $\RR^n$}\label{nil}

Let us consider the euclidean space $\RR^n$ as a differentiable manifold, and let
$\mca_n$ denote the space of all Riemannian metrics on $\RR^n$.  The tangent space
$\tang_0\RR^n$ at the point $0\in\RR^n$ is naturally identified with $\RR^n$, and if
$L(x):\RR^n\longrightarrow\RR^n$ is the translation by $x\in\RR^n$ (i.e. $L(x)y=x+y$
for all $y\in\RR^n$), then $\dif L(x)|_0:\tang_0\RR^n\longrightarrow\tang_x\RR^n$ is
an isomorphism of vector spaces.  The canonical inner product $\ip$ on $\RR^n$ therefore determines a distinguished
element $g_0\in\mca_n$ given by
\begin{equation}\label{defg0}
g_0(x)(\dif L(x)|_0v,\dif L(x)|_0w):=\la v,w\ra, \qquad\forall
v,w\in\RR^n\equiv\tang_0\RR^n, \quad x\in\RR^n,
\end{equation}
that is, the flat metric.  Any inner product on
$\RR^n$ actually defines in the same way a flat element in $\mca_n$.

Let $\mu$  be a bilinear map
$$
\mu:\RR^n\times\RR^n\longrightarrow\RR^n,
$$
and assume that
$\mu$ is skew symmetric, satisfies the {\it Jacobi identity} (i.e. left
multiplication maps $\ad_{\mu}{x}$ are all derivations of the algebra $(\RR^n,\mu)$)
and that $\mu$ is {\it nilpotent} (i.e. the operators $\ad_{\mu}{x}$ are all
nilpotent). In other words, $(\RR^n,\mu)$ is a nilpotent Lie algebra.  It is
well-known that the simply connected (always connected) nilpotent Lie group $N_{\mu}$
with Lie algebra $(\RR^n,\mu)$ is diffeomorphic to $\RR^n$; moreover, the
exponential map $\exp_{\mu}:\RR^n\longrightarrow N_{\mu}$ is a diffeomorphism. One
can therefore identify $N_{\mu}$ with $\RR^n$ via $\exp_{\mu}$ and use the
Baker-Campbell-Hausdorff formula
$$
\exp_{\mu}(x)\exp_{\mu}(y) = \exp_{\mu}(x+y+p_{\mu}(x,y)), \qquad\forall
x,y\in\RR^n,
$$
where $p_{\mu}:\RR^n\times\RR^n\longrightarrow\RR^n$ is a polynomial function (see for instance
\cite[Section 2.15]{Vrd}), to define a group operation on $\RR^n$ by
$$
x\cdot_{\mu}y:=x+y+p_{\mu}(x,y), \qquad\forall x,y\in\RR^n.
$$
In this way, $(\RR^n,\cdot_{\mu})$ becomes the simply connected
nilpotent Lie group with Lie algebra $(\RR^n,\mu)$.

\begin{remark} If
$x=(x_1,...,x_n)$ and $y=(y_1,...,y_n)$, then the product $x\cdot_{\mu}y$ is
polynomial in $x_1,...,y_n$, and conversely, it is proved in \cite{Lzr} that any
group operation on $\RR^n$ which is polynomial in the coordinates is necessarily a
nilpotent Lie group structure on $\RR^n$.
\end{remark}

It is easy to see that $p_{\mu}(ax,bx)=0$ for all $a,b\in\RR$, $x\in\RR^n$, which gives that $0$ is the identity element of the group $(\RR^n,\cdot_{\mu})$ and that $-x$ is the inverse of any $x\in\RR^n$.

We may define a Riemannian metric attached to each nilpotent Lie bracket $\mu$, analogously to what we did in (\ref{defg0}) for the flat metric $g_0$, which will correspond to $\mu=0$.  If $L_{\mu}(x):\RR^n\longrightarrow\RR^n$ denotes left multiplication by $x$ in the group $(\RR^n,\cdot_{\mu})$ (i.e. $L_{\mu}(x)y=x+y+p_{\mu}(x,y)$ for all $y\in\RR^n$), then $L_{\mu}(x)\in\Diff(\RR^n)$ (with inverse $L_{\mu}(-x)$), and thus $\mu$ defines a Riemannian metric $g_\mu$ on $\RR^n$ by
\begin{equation}\label{defgmu}
g_\mu(x)(\dif L_\mu(x)|_0v,\dif L_\mu(x)|_0w):=\la v,w\ra, \qquad\forall
v,w\in\RR^n\equiv\tang_0\RR^n, \quad x\in\RR^n,
\end{equation}
where $\ip=g_\mu(0)$ is the canonical inner product on $\RR^n$.  Thus $L_\mu(x)$ is automatically an isometry of $g_\mu$ for any $x\in\RR^n$, and so the nilpotent Lie group $(\RR^n,\cdot_\mu)$ acts transitively by isometries on $(\RR^n,g_\mu)$.

The {\it degree of nilpotency} $1\leq k_\mu\leq n-1$ of $\mu$ is defined by
$$
k_\mu:=\min\{
j\in\NN:(\ad_{\mu}x)^j=0,\;\forall x\in\RR^n\},
$$
which is usually referred to, by saying that $\mu$ is {\it $k_\mu$-step nilpotent}.  An explicit formula for $p_{\mu}$ is really hard to get when $k_\mu$ is large, although the low order terms are well known (see \cite[(2.15.17)]{Vrd}):
\begin{align}
p_{\mu}(x,y) =& \unm\mu(x,y)+\tfrac{1}{12}\mu(x,\mu(x,y))-\tfrac{1}{12}\mu(y,\mu(x,y)) \notag \\
& -\tfrac{1}{48}\mu(y,\mu(x,\mu(x,y))) -\tfrac{1}{48}\mu(x,\mu(y,\mu(x,y))) \label{pmu}\\
& + \;\mbox{commutators in five or more terms}.\notag
\end{align}

In any case, we have that
\begin{equation}\label{pmuc}
p_\mu(x,y)=\left( p_\mu^1(x,y),...,p_\mu^n(x,y)\right),
\end{equation}
where $p_\mu^i:\RR^n\times\RR^n\longrightarrow\RR$ are polynomial functions on $(x,y)$ of degree $\leq k_\mu$, whose coefficients are universal polynomial expressions on $\mu$ depending only on $n$, of degree $\leq n-2$ (and actually $\leq k_\mu-1$ for each $\mu$).  By `polynomial on $\mu$' we will always mean polynomial on the coordinates $\mu_{ij}^k$'s of $\mu$ defined by
$$
\mu(e_i,e_j)=\sum_k\mu_{ij}^ke_k.
$$
For each $1\leq m\leq n$ we have that
\begin{align}
p_{\mu}^m(x,y) =& \unm\sum\mu_{ij}^m\, x_iy_j +\tfrac{1}{12}\sum\mu_{ij}^k\mu_{lk}^m \, x_ix_ly_j  -\tfrac{1}{12}\sum\mu_{ij}^k\mu_{lk}^m \, x_iy_jy_l \notag \\
& + \;\mbox{monomials of degree}\,\geq 4.\notag
\end{align}

It follows from (\ref{pmu}) that $p_{\mu}(\vp x,\vp y)=\vp p_{\mu}(x,y)$ for all $x,y\in\RR^n$ and $\vp\in\Aut(\RR^n,\mu)$, and hence $\Aut(\RR^n,\mu)\subset\Aut(\RR^n,\cdot_\mu)$.  Conversely, if $\vp\in\Aut(\RR^n,\cdot_\mu)$, then $\vp$ is linear as $\exp:(\RR^n,\mu)\longrightarrow (\RR^n,\cdot_\mu)$ is the identity map and therefore $\vp=\dif\vp|_0$.  By taking the second derivative at $t=0$ of $p_{\mu}(tx,ty)$ one can easily show that $\vp\in\Aut(\RR^n,\mu)$, and, in consequence,
\begin{equation}\label{aut}
\Aut(\RR^n,\cdot_\mu)=\Aut(\RR^n,\mu).
\end{equation}

We note that $L_\mu(x)$ and its inverse $L_\mu(-x)$ are both polynomial diffeomorphisms of $\RR^n$. When $k_\mu=2$, it follows from (\ref{pmu}) that the product is simply given by $x\cdot_\mu y=x+y+\unm\mu(x,y)$, and hence $L_\mu(x)$ is actually an affine map of $\RR^n$: $L_\mu(x)y=(I+\unm\ad_\mu{x})y+x$ for all $y\in\RR^n$.

It will be useful to have an expression of the metric $g_\mu$ in terms of the canonical global chart $(x_1,...,x_n)$ of $\RR^n$.

\begin{lemma}\label{gmuc}
The coordinates $(g_\mu)_{ij}:\RR^n\longrightarrow\RR$ of the metric $g_\mu$ are given by
$$
(g_\mu)_{ij}(x) = \delta_{ij} + \frac{\partial p_\mu^i}{\partial y_j}(-x,x) +\frac{\partial p_\mu^j}{\partial y_i}(-x,x) + \sum_k\frac{\partial p_\mu^k}{\partial y_i}(-x,x) \frac{\partial p_\mu^k}{\partial y_j}(-x,x).
$$
\end{lemma}

\begin{proof}
By (\ref{defgmu}) we have
$$
(g_\mu)_{ij}(x)= \la\dif L_\mu(x)|_0^{-1}e_i,\dif L_\mu(x)|_0^{-1}e_j\ra = \la\dif L_\mu(-x)|_xe_i,\dif L_\mu(-x)|_xe_j\ra,
$$
and since from (\ref{pmuc}) we get that
\begin{align}
\dif L_\mu(-x)|_xe_i &=\ddt|_0L_\mu(-x)(x+te_i) =\ddt|_0 te_i+p_\mu(-x,x+te_i) \\
& =e_i+\left(\frac{\partial p_\mu^1}{\partial y_i}(-x,x),...,\frac{\partial p_\mu^n}{\partial y_i}(-x,x)\right),\notag
\end{align}
the formula follows.
\end{proof}

For a multiindex $\alpha=(\alpha_1,...,\alpha_n)$ we denote $x^\alpha:=x_1^{\alpha_1}...x_n^{\alpha_n}$ for $x=(x_1,...,x_n)$, and $|\alpha|:=\alpha_1+\dots\alpha_n$.

\begin{corollary}\label{gmucc}
$(g_\mu)_{ij}$ is a polynomial on $x$,
$$
(g_\mu)_{ij}(x)=\sum_\alpha a_{\alpha}^{ij}(\mu)x^\alpha,
$$
of degree $\leq 2(k_\mu-1)$, and each coefficient $a_{\alpha}^{ij}$ is a universal polynomial expression on $\mu$ depending only on $i$, $j$, $\alpha$ and $n$, of degree $|\alpha|\leq 2(n-2)$ (and actually $\leq 2(k_\mu-1)$ for each $\mu$).
\end{corollary}

\begin{example}\label{2pasos}
It is easy to see that if $k_\mu=2$ then
$$
(g_\mu)_{ij}(x)=\delta_{ij}-\unm\sum_k(\mu_{kj}^i+\mu_{ki}^j)x_k +\unc\sum_{kl}\left(\sum_r\mu_{ki}^r\mu_{lj}^r\right)\, x_kx_l.
$$
\end{example}

We use the following notion of convergence for metrics on $\RR^n$ (see for instance \cite[Chapter 3]{libro}).

\begin{definition}\label{conv}
Let $\{ g_k\}_{k\in\NN}$, $g$ be Riemannian metrics on $\RR^n$.  We say that {\it $g_k$ converges in $C^\infty$ to $g$ uniformly on compact sets in $\RR^n$} ($g_k\to g$ for short) if for any compact set $K\subset\RR^n$, $p>0$ and $\epsilon>0$, there exists $k_0=k_0(K,p,\epsilon)$ such that for $k\geq k_0$,
$$
\sup_{1\leq q\leq p}\;\sup_{x\in K} \|\nabla^q(g_k-g)\|_g<\epsilon,
$$
where $\nabla$ is the Levi-Civita connection of $g$ and $\|\cdot\|_g$ denotes the norm in the spaces of sections of the corresponding tensor bundles over $\RR^n$.
\end{definition}

\begin{remark}\label{convr}
By using global coordinates in $\RR^n$, convergence $g_k\to g$ can be rephrased as follows: for any multiindex $\alpha$ the partial derivatives $\partial^\alpha(g_k)_{ij}$ converge to $\partial^\alpha g_{ij}$ uniformly on compact sets of $\RR^n$, as $k\to\infty$.
\end{remark}

\begin{remark}
As usual, for a continuous family of metrics $\{ g_t\}$, convergence $g_t\to g$ as $t\to\infty$ will mean $g_{t_k}\to g$ for any sequence $t_k\to\infty$, or equivalently, we may replace $k_0$ by $t_0=t_0(K,p,\epsilon)$ above and require the condition for $t_0\leq t$.
\end{remark}

\begin{proposition}\label{convmu}
$\mu_k\to\lambda$ in $\nca_n\subset V_n$ (usual vector space topology) if and only if $g_{\mu_k}$ converges in $C^\infty$ to $g_\lambda$ uniformly on compact sets in $\RR^n$.
\end{proposition}

\begin{remark}
In \cite{spacehm}, the relationship between the `algebraic' convergence of brackets and other well-known notions, including local, infinitesimal and pointed or Cheeger-Gromov convergence, is studied in the general homogeneous case.
\end{remark}

\begin{proof}
The coordinates $(g_\mu)_{ij}$ of a metric $g_\mu$ have been described in Corollary \ref{gmucc}.  Since the coefficients $a_\alpha(\mu)$ depend polynomially on $\mu$, we see at once that $g_{\mu_k}\to g_\lambda$ follows from $\mu_k\to\lambda$ by using Remark \ref{convr}.

For the converse assertion, we first note that if $\nabla^\mu$ denotes the Levi-Civita connection of $g_\mu$, then
$$
g_\mu(0)\left((\nabla^\mu_{e_r}e_j)_0,e_i\right) = \unm(\mu_{rj}^i+\mu_{ri}^j+\mu_{ji}^r)
$$
(see for instance \cite[7.27]{Bss}), and if $\alpha$ is the multiindex with $1$ at entry $r$ and $0$ elsewhere, then it is easy to see by using Lemma \ref{gmuc} that
$$
\partial^\alpha(g_\mu)_{ij}(0) = -\unm(\mu_{rj}^i+\mu_{ri}^j)
$$
(see Example \ref{2pasos}).  Therefore, the convergence $g_{\mu_k}\to g_\lambda$ implies that $(\mu_k)_{rj}^i+(\mu_k)_{ri}^j+(\mu_k)_{ji}^r \to\lambda_{rj}^i+\lambda_{ri}^j+\lambda_{ji}^r$ and $(\mu_k)_{rj}^i+(\mu_k)_{ri}^j\to\lambda_{rj}^i+\lambda_{ri}^j$, which gives uniform convergence $(\mu_k)_{ji}^r\to\lambda_{ji}^r$ for all $i,j,r$, as $k\to\infty$, and so $\mu_k\to\lambda$.
\end{proof}

\section{Some technical background}

The space of all skew-symmetric algebras of dimension $n$ is
parameterized by the vector space
$$
V_n:=\lamn=\{\mu:\RR^n\times\RR^n\longrightarrow\RR^n : \mu\; \mbox{bilinear and
skew-symmetric}\}.
$$
Then
$$
\nca_n:=\{\mu\in V_n:\mu\;\mbox{satisfies the Jacobi identity and is nilpotent}\}
$$
is an algebraic subset of $V_n$ as all the required conditions  can be written as zeroes of
polynomial functions.  $\nca_n$ is often called the {\it variety of nilpotent Lie algebras} (of
dimension $n$).  There is a natural linear action of $\G$ on $V_n$ given by
\begin{equation}\label{action}
g.\mu(x,y)=g\mu(g^{-1}x,g^{-1}y), \qquad x,y\in\RR^n, \quad g\in\G,\quad \mu\in V_n.
\end{equation}
It is easily seen that $\nca_n$ is $\G$-invariant, the Lie algebra isomorphism classes are
precisely the $\G$-orbits and the isotropy subgroup $\G_\mu$ equals $\Aut(\RR^n,\mu)$ for any $\mu\in\nca_n$.

The representation $\pi:\g\longrightarrow\End(V_n)$ obtained by differentiation of
(\ref{action}) is given by
\begin{equation}\label{actiong}
\pi(\alpha)\mu =-\delta_{\mu}(\alpha)=\ddt|_{t=0}
e^{t\alpha}.\mu =\alpha\mu(\cdot,\cdot)-\mu(\alpha\cdot,\cdot)-\mu(\cdot,\alpha\cdot),
\end{equation}
for all $\alpha\in\g$, $\mu\in V_n$.  We note that $\delta_{\mu}:\glg_n(\RR)\longrightarrow V_n$ is linear and
$\delta_{\mu}(\alpha)=0$ if and only if $\alpha\in\Der(\mu)$, the Lie algebra of
derivations of the algebra $\mu$.

The canonical inner product $\ip$ on $\RR^n$
determines an $\Or(n)$-invariant inner product on $V_n$, also denoted by $\ip$, as
follows:
\begin{equation}\label{innV}
\la\mu,\lambda\ra= \sum\la\mu(e_i,e_j),\lambda(e_i,e_j)\ra
=\sum\la\mu(e_i,e_j),e_k\ra\la\lambda(e_i,e_j),e_k\ra,
\end{equation}
and also the standard $\Ad(\Or(n))$-invariant inner product on $\g$ given by
$$
\la \alpha,\beta\ra=\tr{\alpha \beta^{\mathrm t}}=\sum\la\alpha e_i,\beta e_i\ra
=\sum\la\alpha e_i,e_j\ra\la\beta e_i,e_j\ra, \qquad \alpha,\beta\in\g,
$$
where $\{ e_1,...,e_n\}$ denotes the canonical basis of $\RR^n$.

\begin{remark} We have made
several abuses of notation concerning inner products.  Recall that $\ip$ has been
used in this paper to denote an inner product on $\RR^n$, $V_n$ and $\g$.
\end{remark}

We note that
$\pi(\alpha)^t=\pi(\alpha^t)$ and $(\ad{\alpha})^t=\ad{\alpha^t}$ for any
$\alpha\in\g$, due to the choice of the canonical inner products everywhere.

\section{Geometry of the metrics $g_{\mu}$}\label{geom}

We describe in this section some well-known facts on the geometry of a metric $g_\mu\in\mca_n$, concerning mainly isometries and Ricci curvature.

For each $\mu\in\nca_n$, let $g_\mu$ be the metric on $\RR^n$ defined as in (\ref{defgmu}).  Any left invariant metric on any $n$-dimensional simply connected nilpotent Lie group is isometric to at least one element in the subset
$$
\{g_\mu:\mu\in\nca_n\}\subset\mca_n.
$$
If $\mu,\lambda\in\nca_n$, then $g_\mu$ and $g_\lambda$ are respectively isometric to two left invariant metrics on the same nilpotent Lie group if and only if $\lambda\in\G.\mu$.  For any $\mu\in\nca_n$ and inner product $\ipp$ on $\RR^n$ we may also define $g_{\mu,\ipp}$ as in {\rm (\ref{defgmu})} by using $\ipp$ instead of $\ip$.  The set
$$
\{ g_{\mu,\ipp}:\mu\in\nca_n, \;\ipp \;\mbox{inner product on}\; \RR^n\}\subset\mca_n
$$
is therefore the set of all metrics on $\RR^n$ which are invariant by some transitive nilpotent Lie group.

\begin{theorem}\label{nilprop}\cite{Wls}
\begin{itemize}
\item[(i)] The isometry group is given by $\Iso(\RR^n,g_{\mu})=K_\mu\ltimes L_\mu(\RR^n)$ for any $\mu\in\nca_n$, where $K_\mu$ is the isotropy subgroup at $0\in\RR^n$, and $K_\mu:=\Aut(\RR^n,\mu)\cap\Or(n)$ is the group of $\ip$-orthogonal automorphisms (recall {\rm (\ref{aut})}).

\item[(ii)] If $\mu,\lambda\in\nca_n$, then $g_\mu$ and $g_\lambda$ are isometric if and only if $\lambda\in\Or(n).\mu$ (in particular, they must be invariant under isomorphic nilpotent Lie groups).

\item[(iii)] $g_{\mu,\ipp}$ and $g_{\mu',\ipp'}$ are isometric if and only if there exists $h\in\G$ such that $\mu'=h.\mu$ and $\ipp'=(h\cdot,h\cdot)$.
\end{itemize}
\end{theorem}

It follows from (\ref{defgmu}) that not only the metric $g_{\mu}$ is completely determined by its value at $0$ but also are its scalar curvature
$$
\scalar_\mu:=\scalar(g_{\mu})(0)\in\RR,
$$
its Ricci tensor
$$
\ricci_\mu:=\ricci(g_{\mu})(0):\RR^n\times\RR^n\longrightarrow\RR \qquad \mbox{(symmetric form)}
$$
and its curvature tensor
$$
\Riem_\mu:=\Riem(g_{\mu})(0):\RR^n\times\RR^n\times\RR^n\times\RR^n \longrightarrow\RR \qquad \mbox{($4$-linear map)}.
$$
We denote by
$$
\Ricci_\mu:=\Ricci(g_{\mu})(0):\RR^n\longrightarrow\RR^n \qquad \mbox{(symmetric linear operator)}
$$
the Ricci operator of $g_{\mu}$, which is given by $\ricci_\mu(x,y)=\la\Ricci_\mu x,y\ra$ for all $x,y\in\RR^n$.

\begin{lemma}\label{nilprop2}
Let $\mu\in\nca_n$.
\begin{itemize}
\item[(i)] The Ricci tensor and Ricci operator of $g_\mu$ are respectively given by
$$
\ricci_\mu(x,y)=-\unm\sum\la\mu(x,e_i),e_j\ra\la\mu(y,e_i),e_j\ra +\unc\sum\la\mu(e_i,e_j),x\ra\la\mu(e_i,e_j),y\ra,
$$
for all $x,y\in\RR^n$, and
$$
\Ricci_\mu=-\unm\sum(\ad_{\mu}{e_i})^t\ad_{\mu}{e_i} +
\unc\sum\ad_{\mu}{e_i}(\ad_{\mu}{e_i})^t,
$$
where $\{ e_i\}$ is the canonical basis of $\RR^n$.

\item[(ii)] The scalar curvature of $g_\mu$ is given by $\scalar_\mu=-\unc\|\mu\|^2$ (see {\rm(\ref{innV})}).

\item[(iii)] \cite{Mln} If $\mu\ne 0$ then the Ricci tensor $\ricci_\mu$ has always both positive and negative directions.
\end{itemize}
\end{lemma}

\begin{proof}
Parts (i) and (ii) follow easily from the general curvature formulas for homogeneous spaces (see for instance \cite[7.38]{Bss}).  We prove (iii).  It follows from (i) that any $x\in\RR^n$ orthogonal to $\mu(\RR^n,\RR^n)$ and such that $\ad_\mu{x}\ne 0$ gives $\ricci(x,x)<0$, and any $y\in\mu(\RR^n,\RR^n)$ with $\ad_\mu{y}=0$ has $\ricci(y,y)>0$  (both such $x$ and $y$ are easily seen to exist by the nilpotency of $\mu$).
\end{proof}

\section{Ricci flow starting at a metric $g_{\mu}$}\label{trick}

For a given $\mu_0\in\nca_n$, let $g(t)$ be a {\it Ricci flow} starting at $g_{\mu_0}$, that is, a one-parameter family $\{ g(t)\}\subset\mca_n$ that satisfies the Ricci flow equation
\begin{equation}\label{RF}
\dpar g(t)=-2\ricci(g(t)),\qquad g(0)=g_{\mu_0}.
\end{equation}
The short time existence of a solution follows from \cite{Shi}, as $g_{\mu_0}$ is homogeneous and hence complete and of bounded curvature. Alternatively, one may require $(\RR^n,\cdot_{\mu_0})$-invariance of $g(t)$ for all $t$, and thus $g(t)$ will also be determined by its value at $0$, denoted by $\ip_t:=g(t)(0)$, as in (\ref{defgmu}):
$$
g(t)(x)(\dif L_{\mu_0}(x)|_0v,\dif L_{\mu_0}(x)|_0w):=\la v,w\ra_t, \qquad\forall
v,w\in\RR^n\equiv\tang_0\RR^n, \quad x\in\RR^n.
$$
The Ricci flow equation (\ref{RF}) is therefore equivalent to the ODE
\begin{equation}\label{RFip}
\ddt \ip_t=-2\ricci(\ip_t), \qquad\ip_0=\ip,
\end{equation}
where $\ricci(\ip_t):=\ricci(g(t))(0)$, and hence short time existence and uniqueness of the solution in the class of $(\RR^n,\cdot_{\mu_0})$-invariant metrics is guaranteed.  In this way, $g(t)$ is homogeneous and so complete and of bounded curvature for all $t$, hence the uniqueness also follows from \cite{ChnZhu}.  It is actually a simple matter to prove that the uniqueness, in turn, implies our assumption of $(\RR^n,\cdot_{\mu_0})$-invariance.  The need for this circular argument is due to the fact that the uniqueness of the Ricci flow solution is still an open problem in the noncompact general case (see \cite{Chn}).

On the other hand, for a given $\mu_0\in\nca_n$, we may consider the {\it bracket flow} defined for a curve $\{\mu=\mu(t)\}\subset V_n$ by the following ODE:
\begin{equation}\label{BF}
\ddt\mu=\delta_{\mu}(\Ricci_{\mu}), \qquad \mu(0)=\mu_0,
\end{equation}
where $\delta_\mu:\g\longrightarrow V_n$ is defined by
$$
\delta_{\mu}(\alpha) =\mu(\alpha\cdot,\cdot)+\mu(\cdot,\alpha\cdot)-\alpha\mu(\cdot,\cdot),
\qquad \alpha\in\g,\quad\mu\in V_n.
$$
Equation (\ref{BF}) is well defined as $\Ricci_\mu$ can be computed for any $\mu\in V_n$ as in Lemma \ref{nilprop2}, (i), and not only for $\mu\in\nca_n$.  However, this technicality is only needed to define the ODE, since the solution $\mu(t)$ stays in $\nca_n$ (and even in the submanifold $\G.\mu_0$), as long as $\mu_0\in\nca_n$.  Indeed,
$$
\delta_{\mu}(\Ricci_\mu)=\ddt\mid_{t=0}\!
e^{-t\Ricci_\mu}.\mu\in\tang_\mu\G.\mu, \qquad\forall\mu\in V_n,
$$
and thus $\mu(t)\in\G.\mu_0\subset\nca_n$ for all $t$ by a standard ODE theory argument.

We now show that the flows (\ref{RF}) and (\ref{BF}) are intimately related.  Let $g_{\mu_0}$ be the Riemannian metric on $\RR^n$ defined as in {\rm (\ref{defgmu})} for a nilpotent Lie bracket $\mu_0:\RR^n\times\RR^n\longrightarrow\RR^n$.

\begin{theorem}\label{eqfl}
If $g(t)$ is the solution to the Ricci flow with $g(0)=g_{\mu_0}$ (see {\rm (\ref{RF})}), and $\mu(t)$ is the solution to the bracket flow with $\mu(0)=\mu_0$ (see {\rm (\ref{BF})}), then there exists a one-parameter family $\{ h(t)\}\subset\G$ such that
$$
g(t)=h(t)^*g_{\mu(t)}, \qquad\forall t.
$$
Furthermore, the following conditions hold for $h=h(t)$ for all $t$:

\begin{itemize}
\item[(i)] $\ddt h=-h\Ricci(g(t))(0)=-\Ricci_{\mu(t)}h$, $\quad h(0)=I$.

\item[(ii)] $\ip_t=\la h\cdot,h\cdot\ra$ is the solution to {\rm (\ref{RFip})}.

\item[(iii)] $\mu(t)=h\mu_0(h^{-1}\cdot,h^{-1}\cdot)$.
\end{itemize}
\end{theorem}

\begin{remark}
 The Ricci flow $g(t)$ and the bracket flow $\mu(t)$ differ therefore only by pullback by a time-dependent diffeomorphism (linear map), and are equivalent in the following sense: each one can be obtained from the other by solving the corresponding ODE in part (i) and applying parts (ii) or (iii), accordingly.  In \cite{homRF}, this result has been generalized for Ricci flows starting at any homogeneous Riemannian manifold.
\end{remark}

\begin{remark}\label{GP}
Part (iii) of this theorem has been proved in \cite[Section 2.2]{Gzh}, and we also refer to \cite{Pyn} for a study of the Ricci flow on nilmanifolds via the bracket flow (\ref{BF}).  In the case of $3$-dimensional unimodular Lie groups, a global picture of the qualitative behavior of the Ricci flow is given in \cite{GlcPyn} by using the same approach proposed in the above theorem: to vary brackets instead of inner products. Every metric is represented \`a la Milnor by a triple $(a_1,a_2,a_3)\in\RR^3$ of Lie bracket structural constants in \cite[Theorem 2.4]{GlcPyn}, and the ODE for the normalized quantities $(a_2/a_1,a_3/a_1)$ which is equivalent to the Ricci flow in the sense above is given in \cite[Proposition 3.1]{GlcPyn}.
\end{remark}

\begin{remark}\label{his}
It follows that the map $h(t)\in\G$ satisfies, for all $t$, the following properties:
\begin{itemize}
\item $h(t):(\RR^n,g(t))\longrightarrow(\RR^n,g_{\mu(t)})$ is an isometry.

\item $h(t):(\RR^n,\mu_0)\longrightarrow(\RR^n,\mu)$ is a Lie algebra isomorphism.

\item $h(t):(\RR^n,\cdot_{\mu_0})\longrightarrow(\RR^n,\cdot_\mu)$ is a Lie group isomorphism.

\item $h(t):(\RR^n,\ip_t)\longrightarrow(\RR^n,\ip)$ is an isometry of inner product vector spaces.
\end{itemize}
\end{remark}

\begin{proof}
We first consider the solution $h=h(t)\in\G$ to the ODE
$$
\ddt h=-h\Ricci(\ip_t),\qquad h(0)=I,
$$
where $\Ricci(\ip_t):=\Ricci(g(t))(0)$, which is defined on the same interval of time as $g(t)$ by a standard result in ODE theory ($h(t)$ is easily seen to be invertible for all $t$).  If $\ipp_t:=\la h(t)\cdot,h(t)\cdot\ra$ and $h':=\ddt h(t)$ then
\begin{align*}
\ddt\ipp_t &= \la h'\cdot,h\cdot\ra+\la h\cdot,h'\cdot\ra \\
&= -\la h\Ricci(\ip_t)\cdot,h\cdot\ra -\la h\cdot,h\Ricci(\ip_t)\cdot\ra \\
&= -(\Ricci(\ip_t)\cdot,\cdot)_t -(\cdot,\Ricci(\ip_t)\cdot)_t.
\end{align*}
On the other hand, since $\Ricci(\ip_t)$ is symmetric with respect to $\ip_t$, it follows from (\ref{RFip}) that $\ip_t$ satisfies
$$
\ddt\ip_t=-2\ricci(\ip_t)=-2\la\Ricci(\ip_t)\cdot,\cdot\ra_t= -\la\Ricci(\ip_t)\cdot,\cdot\ra_t-\la\cdot,\Ricci(\ip_t)\cdot\ra_t.
$$
Thus $\ipp_t$ and $\ip_t$, as curves in the manifold $\G/\Or(n)$ of inner products on $\RR^n$, satisfy the same ODE and $\ipp_0=\ip_0=\ip$.  Part (ii) therefore holds by uniqueness of the solution.

It now follows from (ii) that $h(t):(\RR^n,g(t))\longrightarrow(\RR^n,g_{\lambda(t)})$ is an isometry for the curve $\lambda(t):=h(t).\mu_0$ (see Theorem \ref{nilprop}, (iii)).  This implies that $\Ricci_{\lambda(t)}=h(t)\Ricci(\ip_t)h(t)^{-1}$, or equivalently, $h'=-\Ricci_{\lambda(t)}h$, and thus
\begin{align*}
\ddt\lambda &= h'\mu_0(h^{-1}\cdot,h^{-1}\cdot) -h\mu_0(h^{-1}h'h^{-1}\cdot,h^{-1}\cdot) -h\mu_0(h^{-1}\cdot,h^{-1}h'h^{-1}\cdot) \\
&= (h'h^{-1})h\mu_0(h^{-1}\cdot,h^{-1}\cdot)-h\mu_0(h^{-1}(h'h^{-1})\cdot,h^{-1}\cdot) -h\mu_0(h^{-1}\cdot,h^{-1}(h'h^{-1})\cdot) \\
&= -\delta_{\lambda}(h'h^{-1}) =\delta_{\lambda}(\Ricci_\lambda).
\end{align*}
But $\ddt\mu=\delta_\mu(\Ricci_\mu)$ and $\mu(0)=\lambda(0)=\mu_0$, so that $\mu(t)=\lambda(t)$ for all $t$, from which parts (i) and (iii) follow.  We also obtain that $h(t):(\RR^n,g(t))\longrightarrow(\RR^n,g_{\mu(t)})$ is an isometry, concluding the proof of the theorem.
\end{proof}

\section{The bracket flow}

Let $g_{\mu_0}$ be the Riemannian metric on $\RR^n$ corresponding to the nilpotent Lie bracket $\mu_0\in\nca_n$ (see (\ref{defgmu})).  According to Theorem \ref{eqfl}, the Ricci flow $g(t)$ with $g(0)=g_{\mu_0}$ is
equivalent in a natural and specific sense to the ODE (\ref{BF}) for a curve $\mu(t)\in V_n$.  In particular, the maximal interval of
time where a solution exists is the same for both flows.  We also note that at each
time $t$, the Riemannian manifolds $(\RR^n,g(t))$ and $(\RR^n,g_{\mu(t)})$ are isometric, so that the behavior
of the curvature and of any other Riemannian invariant along the
Ricci flow can be studied on the bracket flow $\mu(t)$ given by (\ref{BF}).

It is proved in \cite[Lemma 4.1]{soliton} that the gradient of the functional
$$
F:V_n\longrightarrow\RR, \qquad F(\mu)=\tr{\Ricci_{\mu}^2},
$$
is given by
\begin{equation}\label{gradF}
\grad(F)_{\mu}=-\delta_{\mu}(\Ricci_{\mu}), \qquad\forall\mu\in\lamngo.
\end{equation}
One therefore obtains that the bracket flow (\ref{BF}) is precisely the negative gradient flow of $F$.  It follows from Theorem \ref{eqfl} that
\begin{quote}
the Ricci flow starting at any simply connected nilmanifold is, up to pull-back by time-dependent diffeomorphisms, the negative gradient flow of the square norm of the Ricci tensor on the set of metrics whose value at the point $0\in\RR^n$ is $g(0)$ and are invariant under some nilpotent Lie group.
\end{quote}
This remarkable fact paves the way to find estimates and get convergence results which may not be evident in the genuine Ricci flow equation.

\begin{remark}
It is proved in \cite{minimal} that the moment map $m:\nca_n\longrightarrow\g$ for the $\G$-action (\ref{action}) (see \cite[(4.4)]{cruzchica}) is given by
$$
m(\mu)=\tfrac{4}{\|\mu\|^2}\Ricci_\mu.
$$
Thus $F$ is, up to scaling, the square norm of the moment map and (\ref{BF}) its negative gradient flow.
\end{remark}

\subsection{Some ODE's along the bracket flow}
We are interested in the behavior of the Ricci and scalar curvature along the bracket flow.

\begin{lemma}\label{prop}
Let $\mu\in\nca_n$.
\begin{itemize}
\item[(i)] $\delta_{\mu}(I)=\mu$.

\item[(ii)] $\delta_{\mu}^t(\mu)=-4\Ricci_{\mu}$, where $\delta_{\mu}^t:V_n\longrightarrow\g$ is the transpose of $\delta_{\mu}$.  Equivalently,
$$
\tr{\Ricci_\mu\alpha}=-\unc\la\delta_\mu(\alpha),\mu\ra, \qquad\forall\alpha\in\g.
$$

\item[(iii)] $\tr{\Ricci_{\mu}D}=0$ for any $D\in\Der(\RR^n,\mu)$.
\end{itemize}
\end{lemma}

\begin{proof}
It follows from \cite[Proposition 3.5]{minimal} and Lemma \ref{nilprop2}, (i) that
$$
\la 2\Ricci_\mu,\alpha\ra=\unm\la\pi(\alpha)\mu,\mu\ra, \qquad\forall\alpha\in\g,
$$
thus (ii) holds by using (\ref{actiong}). Part (iii) follows from (ii) and part (i) does from the definition of $\delta_\mu$.
\end{proof}

Recall that $\delta_\mu(\Ricci_\mu)=0$ if and only if $\Ricci_\mu\in\Der(\RR^n,\mu)$, but according to part (iii) in the above lemma this is equivalent to $\Ricci_\mu=0$.  From Lemma \ref{nilprop2}, (ii), we conclude that the only fixed point of the bracket flow (\ref{BF}) is $\mu=0$, the flat metric.

\begin{lemma}\label{eqs}
The bracket flow {\rm (\ref{BF})} for $\mu=\mu(t)$ produces the following ODE's:
\begin{itemize}
\item[(i)] $\ddt\Ricci_\mu=-\unm\Delta_\mu(\Ricci_\mu)$, where
$$
\Delta_{\mu}:=S\circ\delta_{\mu}^t\delta_{\mu}:\g\longrightarrow\g, \qquad S(\alpha):=\unm(\alpha+\alpha^t).
$$

\item[(ii)] $\ddt\scalar_\mu=2\tr{\Ricci_\mu^2}=2||\ricci_\mu||^2$.

\item[(iii)] $\ddt \|\mu\|^2=-8\tr{\Ricci_\mu^2}$.
\end{itemize}
\end{lemma}

\begin{remark}
The operator $\delta_{\mu}^t\delta_{\mu}$ is precisely the Laplacian on $$C^1(\RR^n,\mu)=(\RR^n)^*\otimes\RR^n=\g$$ of the Lie algebra cohomology of $(\RR^n,\mu)$ relative to the adjoint representation.
\end{remark}

\begin{proof}
Let us also denote by $\Ricci$ the map defined by
$$
\Ricci:V_n\longrightarrow\g, \qquad\Ricci(\mu):=\Ricci_{\mu}.
$$
It follows from Lemma \ref{prop}, (ii) and (\ref{actiong}) that
\begin{align*}
\la\dif\Ricci|_{\mu}\delta_{\mu}(\alpha),\beta\ra &= \la\ddt|_0
\Ricci_{e^{-t\alpha}.\mu},\beta\ra
= \unc\ddt|_0\la\pi(\beta)e^{-t\alpha}.\mu,e^{-t\alpha}.\mu\ra \\
&= \unm\la\pi(\beta)\mu,\delta_{\mu}(\alpha)\ra =
-\unm\la\delta_{\mu}^t\delta_{\mu}(\alpha),\beta\ra, \qquad\forall\beta\in\g,
\quad \beta^t=\beta,
\end{align*}
from where one deduces that
\begin{equation}\label{dR}
\dif\Ricci|_{\mu}\delta_{\mu}(\alpha)=-\unm\Delta_{\mu}(\alpha), \qquad\forall
\alpha\in\g.
\end{equation}
This implies part (i) since
$$
\ddt\Ricci_\mu=\dif\Ricci|_\mu\ddt\mu= \dif\Ricci|_\mu\delta_\mu(\Ricci_\mu)= -\unm\Delta_\mu(\Ricci_\mu).
$$
Part (iii) follows from (ii) by using that $\tr{\Ricci_\mu}=-\unc\|\mu\|^2$ (see Lemma \ref{nilprop2}, (ii)), and by using (i) we can prove part (ii) as follows:
\begin{align*}
\ddt\scalar_\mu &= \ddt\tr{\Ricci_\mu}=-\unm\tr{\Delta_\mu(\Ricci_\mu)} =-\unm\tr{\delta_\mu^t\delta_\mu(\Ricci_\mu)} \\
&= -\unm\la\delta_\mu^t\delta_\mu(\Ricci_\mu),I\ra= -\unm\la\Ricci_\mu,\delta_\mu^t\delta_\mu(I)\ra =2\la\Ricci_\mu,\Ricci_\mu\ra.
\end{align*}
We have used in the last equality Lemma \ref{prop}, (i) and (ii).
\end{proof}

\subsection{Type-III solutions}\label{t3}
By using Theorem \ref{eqfl} and the ODE for $\|\mu\|^2$ obtained above, we are now in a position to prove that the Ricci flows on simply connected nilmanifolds are all type-III solutions.

\begin{definition}\label{deft3}
A Ricci flow $g(t)$ on a manifold is said to be a {\it type-III solution} if
it is defined for $t\in[0,\infty)$ and there exists $C\in\RR$ such that
$$
||\Riem(g(t))||\leq \frac{C}{t}, \qquad\forall t\in(0,\infty).
$$
\end{definition}

It follows from Lemma \ref{eqs}, (iii), that a solution $\mu(t)$ to (\ref{BF}) will stay for ever in a compact subset,
which implies that $\mu(t)$ is defined for all $t\in[0,\infty)$ for any starting
point $\mu_0$.  Furthermore, it follows from Lemma \ref{nilprop2}, (ii) that
$$
\ddt \|\mu\|^2=-8\tr{\Ricci_{\mu}^2}\leq-\tfrac{8}{n}(\tr{\Ricci_{\mu}})^2=-\tfrac{1}{2n}\|\mu\|^4,
$$
which implies
$$
\|\mu\|^2\leq\frac{1}{\tfrac{t}{2n}+\tfrac{1}{\|\mu_0\|^2}}\leq \frac{2n}{t}, \qquad\forall
t,
$$
and hence we get that
\begin{equation}\label{muto0}
\lim_{t\to\infty}\mu(t)=0.
\end{equation}
The Riemann curvature tensor of $g_\mu$ satisfies $\Riem_{c\mu}=c^2\Riem_{\mu}$ for
any $c\in\RR$ (see \cite[7.30]{Bss}).  Thus
$$
\|\Riem_{\mu(t)}\|= \|\mu\|^2\|\Riem_{\mu/\|\mu\|}\|\leq \frac{2nM_n}{t},
$$
where $M_n$ is the maximum of the continuous function
$\lambda\mapsto\|\Riem_{\lambda}\|$ restricted to the unit sphere of $V_n$.

The following result thus follows, as an application of Theorem \ref{eqfl}.

\begin{theorem}\label{rfnil}
For any $\mu_0\in\nca_n$, the Ricci flow $g(t)$ with $g(0)=g_{\mu_0}$ is a type-III solution for a constant $C_n$ that only depends on the dimension $n$.
\end{theorem}

\begin{corollary}
Let $N$ be a   simply connected nilpotent Lie group endowed with a left invariant
metric $g$.  Then the Ricci flow $g(t)$ with $g(0)=g$ is a type-III solution for a constant $C_n$ that only depends on $n=\dim{N}$.
\end{corollary}

For the Ricci tensor, we know that the maximum of $\lambda\mapsto\|\ricci_{\lambda}\|$
restricted to the unit sphere of $V_n$ is $\tfrac{\sqrt{3}}{4}$ (see
\cite[Theorem 4.6]{strata}), and so
$$
\|\Ricci_{\mu(t)}\|= \|\mu(t)\|^2\|\Ricci_{\mu/\|\mu\|}\|\leq \frac{\sqrt{3}n}{2t}, \qquad\forall t\in(0,\infty).
$$

\subsection{Ricci solitons}
A complete Riemannian metric $g$ on a differentiable manifold $M$ is said to be a
{\it Ricci soliton} if its Ricci tensor satisfies
\begin{equation}\label{defrs}
\ricci(g)=cg+\lca_Xg, \qquad\mbox{for some}\quad c\in\RR, \quad
X\in\chi(M)\;\mbox{complete},
\end{equation}
where $\chi(M)$ denotes the space of all differentiable vector fields on $M$ and
$\lca_X$ the Lie derivative.  If in addition
$X$ is the gradient field of a smooth function $f:M\longrightarrow\RR$, then
(\ref{defrs}) becomes $\ricci(g)=cg+2\Hess(f)$ and $g$ is called a {\it gradient}
Ricci soliton. The main significance of the concept is that $g$ is a Ricci soliton if and
only if the curve of metrics
\begin{equation}\label{rssol}
g(t)=(-2ct+1)\vp_t^*g,
\end{equation}
is a solution to the Ricci flow for some one-parameter group $\vp_t$ of diffeomorphisms of $M$.  According to (\ref{rssol}), Ricci solitons are called {\it expanding}, {\it steady}, or {\it shrinking} depending on
whether $c<0$, $c=0$, or $c>0$ (see \cite[Chapter I]{libro} for further information on Ricci solitons).

The only critical point of the functional $F:V_n\longrightarrow\RR$, $F(\mu)=\tr{\Ricci_{\mu}^2}$, is $\mu=0$, as $F$ is a homogeneous polynomial (of degree $4$) vanishing only at $\mu=0$ (see Lemma \ref{nilprop2}, (ii)).  Let us in turn consider the cone
$$
\cca_n:=\{\mu\in V_n:[\mu] \;\mbox{is a critical point of}\; F:\PP V_n\longrightarrow\RR\},
$$
where $[\mu]$ denotes the class of $\mu$ in the projective space $\PP V_n$.  We note that $\mu\in\cca_n$ if and only if $\mu$ is a critical point of $F$ restricted to the sphere that $\mu$ belongs to.

\begin{theorem}\cite{soliton}\label{MA}
The following conditions are equivalent for a metric $g_\mu$, $\mu\in\nca_n$:
\begin{itemize}
\item[(i)] $g_\mu$ is a Ricci soliton.

\item[(ii)] $\mu\in\cca_n$.

\item[(iii)] $\Ricci_\mu=cI+D$ for some $c\in\RR$, $D\in\Der(\RR^n,\mu)$.
\end{itemize}
Furthermore, if $\mu\in\cca_n$, then $\cca_n\cap\G.\mu=\RR^*\Or(n).\mu$.
\end{theorem}

\begin{remark}
From the last assertion we deduce that a nilpotent Lie group can admit at most one Ricci soliton left invariant metric up to isometry and scaling (see Theorem \ref{nilprop}).
\end{remark}

A Ricci soliton $g_\mu$ is often called a {\it nilsoliton} in the literature, and when nonflat, it is always expanding (indeed, $c=-\tfrac{\tr{\Ricci_\mu^2}}{4\|\mu\|^4}<0$ by Theorem \ref{MA}, (iii), Lemma \ref{prop}, (iii) and Lemma \ref{nilprop2}, (ii)) and it can never be gradient (see \cite{PtrWyl}).  Nilsolitons have been extensively studied because of their interplay with Einstein solvmanifolds (see the survey \cite{cruzchica}).  All known nontrivial examples of homogeneous Ricci solitons are isometric to left-invariant metrics on solvable Lie groups satisfying condition (iii) in Theorem \ref{MA}, so called {\it solsolitons}.  It is proved in \cite{solsolitons} that all solsolitons can be obtained as semidirect products of nilsolitons and an abelian group of symmetric automorphisms.

\section{Normalized Ricci flows and convergence}\label{norm}

We have seen in (\ref{muto0}) that the bracket flow $\mu(t)$ defined in (\ref{BF}) converges to $0$, as $t\to\infty$, for any initial condition $\mu_0\in\nca_n$.  It follows from Proposition \ref{convmu} that $g_{\mu(t)}$ converges in $C^\infty$ to the flat metric $g_0$, and so Theorem \ref{eqfl} shows that,
\begin{quote}
up to pull-back by time-dependent diffeomorphisms, the Ricci flow starting at any simply connected nilmanifold converges in $C^\infty$ to the flat metric uniformly on compact sets in $\RR^n$.
\end{quote}
In order to get a more interesting convergence behavior, we are forced to consider appropriate normalizations of the flows.

\subsection{Normalized Ricci flows}
Let $g_{\mu_0}$ be the Riemannian metric on $\RR^n$ corresponding to the nilpotent Lie bracket $\mu_0\in\nca_n$ (see (\ref{defgmu})).  By rescaling the metric and reparametrizing the time
variable $t$, one can transform the Ricci flow (\ref{RF}) into an {\it $r$-normalized Ricci flow}
\begin{equation}\label{RFrn}
\dpar g(t)=-2\ricci(g(t))-2r(t)g(t),\qquad g(0)=g_{\mu_0},
\end{equation}
for some function $r:[0,T)\longrightarrow\RR$ which may depend on $g(t)$.  A scalar Riemannian invariant of the solution $g(t)$ to (\ref{RFrn}) may remains constant in time as a result of an appropriate choice of the function $r(t)$, as this is actually the flow equation the family $c(s)g(s)$ satisfies for some scaling $c(s)>0$, $c(0)=1$, and $g(s)$ the solution to the unnormalized Ricci flow $\frac{\partial}{\partial s} g(s)=-2\ricci(g(s))$ (the time reparametrization is given by $t(s):=\int_0^s c(u)\;du$ and $r(t)=-\unm\frac{c'(s)}{c(s)^2}$).  If $\ip_t:=g(t)(0)$, then the flow (\ref{RFrn}) is equivalent to
\begin{equation}\label{RFiprn}
\ddt\ip_t=-2\ricci(\ip_t)-2r(t)\ip_t, \qquad \ip_0=\ip,
\end{equation}
and if we define the {\it $r$-normalized bracket flow} for $\mu=\mu(t)$ and $r=r(t)$ by
\begin{equation}\label{BFrn}
\ddt\mu=\delta_{\mu}(\Ricci_\mu+rI) =\delta_{\mu}(\Ricci_\mu)+r\mu, \qquad \mu(0)=\mu_0,
\end{equation}
then the following result may be proved in much the same way as Theorem \ref{eqfl}.  Let $g_{\mu_0}$ be the Riemannian metric on $\RR^n$ defined as in {\rm (\ref{defgmu})} for a nilpotent Lie bracket $\mu_0:\RR^n\times\RR^n\longrightarrow\RR^n$.

\begin{theorem}\label{eqflrn}
If $g(t)$ is the solution to the $r$-normalized Ricci flow with $g(0)=g_{\mu_0}$ (see {\rm (\ref{RFrn})}), and $\mu(t)$ is the solution to the $r$-normalized bracket flow with $\mu(0)=\mu_0$ (see {\rm (\ref{BFrn})}), then there exists a one-parameter family $\{ h(t)\}\subset\G$ such that
$$
g(t)=h(t)^*g_{\mu(t)}, \qquad\forall t.
$$
Furthermore, the following conditions hold for $h=h(t)$ for all $t$:

\begin{itemize}
\item[(i)] $\ddt h=-h(\Ricci(g(t))(0)+r(t)I)=-(\Ricci_{\mu(t)}+r(t)I)h$, $\quad h(0)=I$.

\item[(ii)] $\ip_t=\la h\cdot,h\cdot\ra$ is the solution to {\rm (\ref{RFiprn})}.

\item[(iii)] $\mu(t)=h\mu_0(h^{-1}\cdot,h^{-1}\cdot)$.
\end{itemize}
\end{theorem}

It is worth mentioning at this point that by Proposition \ref{convmu}, any convergence $\mu(t)\to\lambda$ we may get for some flow of the form (\ref{BFrn}) gives rise to a convergence $g_{\mu(t)}\to g_\lambda$ in $C^\infty$.  We first show that the possible limits of any of these normalized Ricci flows are all solitons.

\begin{proposition}\label{limrs}
Let $\mu(t)$ be a solution to {\rm (\ref{BFrn})} with maximal interval of time $[0,T)$, and assume that $\mu(t)$ converges as $t\to T$ to an element $\lambda\in V_n$.  Then $T=\infty$, $\lambda\in\nca_n$ and $g_\lambda$ is a Ricci soliton.  If in addition $\lambda\ne 0$ (i.e. $g_\lambda$ nonflat), then the solution $h(t)$ to the ODE in Theorem \ref{eqflrn}, (i), converges exponentially fast to $0$, as $t\to\infty$ .
\end{proposition}

\begin{remark}
The fact that $h(t)\to 0$ is what makes so difficult to visualize where the genuine $r$-normalized Ricci flow $g(t)(0)=\ip_t=\la h(t)\cdot,h(t)\cdot\ra$ is approaching to.
\end{remark}

\begin{proof}
As $\nca_n$ is closed we have $\lambda\in\nca_n$, and since $\mu(t)$ stays in a compact subset of $V_n$, it follows that $T=\infty$.  By assumption, $\ddt\mu(t)\to 0$ as $t\to\infty$, which gives $r(t)\to r_\infty\in\RR$ and $\delta_\lambda(\Ricci_\lambda+r_\infty I)=0$ (i.e.  $\Ricci_\lambda+r_\infty I\in\Der(\RR^n,\lambda)$).  Theorem \ref{MA} now shows that $g_\lambda$ is a Ricci soliton.  Concerning the last assertion, we first note that if $\lambda\ne 0$ then it is known that that the derivation
$\Ricci_\lambda+r_\infty I$ is positive definite (see \cite[Section 4]{Hbr} or \cite[Lemma 2.17]{standard}).  Let $m_\mu$ denote the minimum eigenvalue of $\Ricci_\mu+rI$.  It follows that for sufficiently large $t$, $0<\unm m_\lambda< m_\mu$ and hence
$$
\ddt\| h\|^2=2\la\ddt h,h\ra= -2\la(\Ricci_\mu+rI)h,h\ra\leq -2m_\mu\| h\|^2\leq-m_\lambda\| h\|^2.
$$
This gives $\| h\|^2\leq e^{-m_\lambda t}$, concluding the proof.
\end{proof}

\subsection{Scalar curvature normalization}\label{scn}
In the homogeneous (possibly noncompact) case, a natural geometric quantity to keep unchanged in time along a normalized Ricci flow is the scalar curvature, as it is a single number associated to each metric which therefore does not need to be integrated.  Recall that the scalar curvature of any nonflat metric $g_\mu$ is negative (see Lemma \ref{nilprop2}, (ii)).

If $\|\mu_0\|=2$, i.e. $\scalar(g_{\mu_0})=-1$, then the solution $\mu(t)$ to the $r$-normalized bracket flow (\ref{BFrn}) for $r(t):=\tr{\Ricci_{\mu(t)}^2}$ satisfies $\|\mu(t)\|\equiv 2$.  Indeed, by Lemma \ref{prop}, (ii) we have that
\begin{align*}
\ddt\|\mu\|^2 &= \la\mu',\mu\ra=\la\delta_\mu(\Ricci_\mu)+r\mu,\mu\ra =\la\Ricci_\mu,\delta_\mu^t(\mu)\ra+r\|\mu\|^2 \\
&= -4\tr{\Ricci_\mu^2}+r_\mu\|\mu\|^2 = \tr{\Ricci_\mu^2}(-4+\|\mu\|^2),
\end{align*}
and hence $\|\mu\|^2\equiv 4$ by uniqueness of the solution since $\|\mu_0\|^2=4$.  Thus the scalar curvature satisfies $\scalar(g_{\mu(t)})\equiv -1$.  It follows from Theorem \ref{eqflrn} that the normalized Ricci flow
\begin{equation}\label{RFipn}
\ddt\ip_t=-2\ricci(\ip_t)-2\|\ricci(\ip_t)\|^2\ip_t,
\end{equation}
satisfies $\scalar(\ip_t)\equiv -1$ as soon as $\scalar(\ip_0)=\scalar(g_{\mu_0})=-1$.

\begin{remark}
This normalized Ricci flow $\ip_t$ equals $-\scalar(\ip_s)\ip_s$, where $\ip_s$ denotes the unnormalized Ricci flow and $s$ a suitable reparametrization in time.
\end{remark}

We shall therefore focus, from now on, on the normalized bracket flow
\begin{equation}\label{BFn}
\ddt\mu=\delta_{\mu}(\Ricci_\mu)+\tr{\Ricci_\mu^2}\,\mu,
\end{equation}
which satisfies $\|\mu(t)\|\equiv 2$ if $\|\mu_0\|=2$.  This flow is precisely the negative gradient flow of $F(\mu)=\tr{\Ricci_{\mu}^2}$ restricted to the sphere $S_2:=\{\mu\in V_n:\|\mu\|=2\}$ (see (\ref{gradF})).  By the compactness of $S_2$, $\mu(t)$ is defined for all $t\in[0,\infty)$.  We note that $\cca_n\cap S_2$, is precisely the set of critical points of $F:S_2\longrightarrow\RR$.  On the other hand, since $F:V_n\longrightarrow\RR$ is a homogeneous polynomial such that $F(\mu)>0$ for any nonzero $\mu\in V_n$, it follows that its negative gradient flow starting at any point converges to $0$.  Moreover, it is proved in \cite{Krd} (see also \cite{Mss}), that the radial projection of such a flow on any sphere of $V_n$ has a unique limit. But our solution $\mu(t)$ to (\ref{BFn}) is precisely a reparametrization of such a projection on $S_2$, and therefore $\mu(t)$ converges as $t\to\infty$ to a single critical point $\lambda\in\nca_n\cap\cca_n\cap S_2$, a Ricci soliton.

By applying Proposition \ref{convmu}, we thus obtain the following result.

\begin{theorem}
Let $\mu(t)$ be the solution to {\rm (\ref{BFn})} with $\mu(0)=\mu_0\in\nca_n$,  $\|\mu_0\|=2$.  Then $g_{\mu(t)}$ converges in $C^\infty$ to a Ricci soliton $g_\lambda$ uniformly on compact sets in $\RR^n$, as $t\to\infty$.
\end{theorem}

According to Theorem \ref{eqfl}, we have that
\begin{quote}
 up to rescaling (scalar curvature $\equiv -1$) and pull-back by time-dependent diffeomorphisms, the Ricci flow starting at any simply connected nilmanifold converges in $C^\infty$ to a Ricci soliton metric uniformly on compact sets in $\RR^n$.
\end{quote}

\begin{remark}
It is proved in \cite[Theorem 6.4]{Jbl} that if the nilpotent Lie group $(\RR^n,\cdot_{\mu_0})$ admits a left invariant Ricci soliton (i.e. $\G.\mu_0\cap\cca_n\ne\emptyset$), or equivalently, if there exists an $(\RR^n,\cdot_{\mu_0})$-invariant Ricci soliton $g$ on $\RR^n$, then the limit $\lambda\in\G.\mu_0$, that is, $g_\lambda$ is invariant by a group isomorphic to $(\RR^n,\cdot_{\mu_0})$ and so $g_\lambda$ is isometric to $g$.
\end{remark}

\end{document}